\documentclass[11pt]{article}
\usepackage[T1]{fontenc}
\usepackage{amsmath,amsfonts,amsthm,amssymb,mathtools}
\usepackage{mathrsfs,graphicx,color,latexsym,calc}
\usepackage[colorlinks,bookmarksopen,bookmarksnumbered,citecolor=blue,linkcolor=red,urlcolor=blue]{hyperref}
\usepackage{enumitem}
\usepackage{authblk}
\newcommand{\address}[1]{{\small #1}}
\usepackage{booktabs,array}

\allowdisplaybreaks[2]
\newtheorem{theorem}{Theorem}[section]
\newtheorem{lemma}[theorem]{Lemma}
\newtheorem{proposition}[theorem]{Proposition}
\newtheorem{corollary}[theorem]{Corollary}

\newtheorem{conjecture}{Conjecture}[section]
\newtheorem*{conjecture*}{Conjecture}
\theoremstyle{definition}

\theoremstyle{remark}

\numberwithin{equation}{section}

\newcommand{\one}{\mathbf{1}}
\newcommand{\E}{\mathbb E}

\hypersetup{
  pdftitle={A proof of the
Manickam--Mikl\'os--Singhi conjecture},
  pdfauthor={Liangdong Fan and Yaojun Chen},
  pdfsubject={A proof of the
Manickam--Mikl\'os--Singhi conjecture},
  pdfkeywords={Manickam--Mikl\'os--Singhi conjecture, nonnegative k-sums, stability, sampling without replacement}
}

\begin{document}

\title{A proof of the
Manickam--Mikl\'os--Singhi conjecture}

\author{Liangdong FAN, \, Yaojun CHEN\thanks{\em\protect\raggedright
Corresponding author. Email addresses:
\href{mailto:liangdongfan@smail.nju.edu.cn}{\textcolor{black}{\mbox{liangdongfan@smail.nju.edu.cn}}} (L. Fan),
\href{mailto:yaojunc@nju.edu.cn}{\textcolor{black}{\mbox{yaojunc@nju.edu.cn}}} (Y. Chen).
}\\
\address{School of Mathematics, Nanjing University, Nanjing 210093, P.R. China}\\
}
\date{}
\begingroup
\topsep=\dimexpr\topsep\relax
\partopsep=\dimexpr\partopsep\relax
\parskip=0pt
\maketitle
\vspace{-12.65pt}

\begin{abstract}
Manickam, Mikl\'os and Singhi (1988) conjectured that
whenever $n\ge4k$, every sequence of $n$ real numbers with nonnegative
sum has at least $\binom{n-1}{k-1}$ nonnegative $k$-sums.
In this paper, we prove the conjecture and show the extremal family of nonnegative $k$-sets is unique.
We further establish a quantitative stability theorem.

\end{abstract}
\endgroup

\baselineskip=0.202in

\section{Introduction}\label{sec:intro}
Let $x_1,\ldots,x_n$ be real numbers with nonnegative sum, and fix
$1\le k<n$. How many of their $k$-sums must be nonnegative?
Assume that
$x_1\ge x_2\ge\cdots\ge x_n$.
For $S\subseteq[n]=\{1,\ldots,n\}$, write
$f(S)=\sum_{i\in S}x_i$.
We call $S$ negative if $f(S)<0$ and nonnegative otherwise.
Let $\mathcal G$ be the family of nonnegative $k$-sets, and let
$\mathcal S_1$ be the star consisting of all $k$-sets containing $1$.
For $(x_1,\ldots,x_n)=(n-1,-1,\ldots,-1)$, we have
$\mathcal G=\mathcal S_1$, so
$|\mathcal G|=\binom{n-1}{k-1}$.
This construction motivates the following conjecture~\cite{MM}.
An equivalent formulation uses the first distribution invariant of the
Johnson scheme~\cite{MS}.

\begin{conjecture}[Manickam--Mikl\'os--Singhi, {\cite{MM,MS}}]
\label{conj:mms}
Suppose that $n\ge4k$, and let $x_1,\ldots,x_n$ be real numbers
with $x_1+\cdots+x_n\ge0$. Then at least $\binom{n-1}{k-1}$
of the $k$-element subsets $S\subseteq[n]$ satisfy
$\sum_{i\in S}x_i\ge0$.
\end{conjecture}

The star bound can fail when $n=3k+1$. Take $3k-2$ entries equal to $3$
and the remaining three equal to $-(3k-2)$~\cite{BM,MM}. The total
sum is zero, but a $k$-sum that uses even one of the negative entries
is at most $3(k-1)-(3k-2)=-1$. Thus exactly $\binom{3k-2}{k}$
of the $k$-sums are nonnegative, fewer than $\binom{3k}{k-1}$ for
$k\ge3$. For $3\le k\le7$, this example is extremal at $n=3k+1$;
for $4\le k\le7$, the sharp bound $|\mathcal G|\ge\binom{n-1}{k-1}$
holds for all $n\ge3k+2$~\cite{HS}.
For $k=3$, Manickam~\cite{ManickamThesis,MM} proved the bound for $n\ge11$, and Marino and Chiaselotti~\cite{MC} later gave another proof for $n\ge12$. A later proof also covers $n\ge11$~\cite{Chow};
for a characterization of equality, see~\cite{ChowThesis}.

The conjecture was known to hold when $k$ divides $n$~\cite{BM,MS},
while early results for general $k$ required $n$ to be very large
relative to $k$~\cite{MM,Tyo}.
Alon, Huang and
Sudakov~\cite{AHS} obtained the first polynomial bound, proving the
conjecture for $n\ge33k^2$ by combining fractional covers with
probabilistic estimates. Frankl~\cite{Frankl} gave a short counting
proof for $n\ge3k^3/2$, and Chowdhury, Sarkis and Shahriari~\cite{CSS}
used a matrix argument to prove the conjecture for $n\ge8k^2$.
Both papers also characterized the families attaining equality in their
respective ranges. By averaging over regular hypergraphs,
Pokrovskiy~\cite{Pok} obtained the first linear bound, $n\ge10^{46}k$.

A different line of research concerns the structure of extremal weight
sequences. Examples with only two distinct weights led to the conjecture
that such a sequence always attains the minimum number of nonnegative
$k$-sums~\cite{AB}. For two-valued sequences, the problem reduces to
estimating a hypergeometric tail~\cite{AKP}. In this setting, the star
bound was proved for $n\ge8k$~\cite{AP}; a preprint extends it to
$n\ge4k$~\cite{Anon}.

In this paper, we prove Conjecture \ref{conj:mms} completely. Moreover, we
characterize the extremal family attaining the equality.

\begin{theorem}\label{thm:main}
If $n\ge4k$, then
\(
  |\mathcal G|\ge\binom{n-1}{k-1},
\)
with equality if and only if $\mathcal G=\mathcal S_1$.
\end{theorem}

Alon, Huang and Sudakov~\cite{AHS} proved a Hilton--Milner-type strengthening of the Conjecture \ref{conj:mms}: if \(n>500k^2\) and \(\mathcal G\) contains no full star, then
\(
|\mathcal G|\ge
\binom{n-1}{k-1}+\binom{n-k-1}{k-1}-1.
\)
They also obtained a density bound when every star contains a fixed positive proportion of negative \(k\)-sets. We further establish a quantitative stability theorem showing that near-extremal families of nonnegative \(k\)-sets must be close to \(\mathcal S_1\).

\begin{theorem}\label{thm:stability}
There exists an absolute constant $C>0$ such that,
for all $n\ge4k$,
\[
  |\mathcal G\mathbin\triangle\mathcal S_1|
  \le C\left(|\mathcal G|-\binom{n-1}{k-1}\right).
\]
\end{theorem}

Theorem~\ref{thm:stability} not only implies Theorem~\ref{thm:main}
but also yields the two applications presented in Section~\ref{sec:applications}. The first concerns the maximum number
of edges in a $k$-uniform hypergraph with no perfect fractional matching.
Using fractional matching duality~\cite{AHS,AB,Blin}, we obtain a linear
stability version of the sharp edge bound for all $n\ge4k$. Our second application concerns the probability bound $\Pr(Y\ge0)\ge p$
for centered weighted Bernoulli sums $Y$. Following Cs\'oka~\cite{Csoka},
we apply Theorem~\ref{thm:stability} and let $n\to\infty$ with $k/n\to p$.
This gives a linear stability version of the bound for $0<p\le1/4$.
In particular, the probability bound for finite sums with positive weights
proves the conjecture of Milojevi\'c and Sudakov~\cite{MiloSud} in this range.

\subsection*{Proof strategy}
The main difficulty lies in the range \(4k<n<5k\).  Indeed, the case $k\mid n$
is known~\cite{BM,MS}, and the disjoint-pair argument in~\cite{MM}
gives the star bound on $n$ points once it is known on $n-k$ points.

If every member of the star \(\mathcal S_1\) is nonnegative, the result follows immediately. Let \(b\) denote the number of negative \(k\)-sets in \(\mathcal S_1\), and let \(c\) denote the number of nonnegative \(k\)-sets outside \(\mathcal S_1\). Since
\(
|\mathcal G|=|\mathcal S_1|-b+c,
\)
we may assume that \(b>0\) and \(c\le b\). We shall derive a contradiction.

The counting argument is particularly simple when \(\mathcal S_1\) contains exactly one negative member \(S\). Since \(S\) is negative, its complement has positive total weight and hence contains many nonnegative \(k\)-sets, all lying outside the star. This immediately gives a contradiction. More generally, each negative member of the star yields many nonnegative \(k\)-sets in its complement. The difficulty is that these sets may be counted multiple times, as the same nonnegative \(k\)-set can be disjoint from several negative members of the star.

To control this overlap, we exploit the maximality of \(x_1\). If \(\{1\}\cup B\) is a negative member of the star, then every \(k\)-set containing \(B\) is negative, since replacing \(1\) by any other index cannot increase the sum. Thus, each such \(B\) is disjoint from many nonnegative \(k\)-sets, but is contained in none. Combining these two constraints, we use spectral estimates to show in Section~\ref{sec:complement} that a substantial proportion of the $k$-sets in $\mathcal S_1$ must have negative sums under our assumption.

The negative \(k\)-sets in the star \(\mathcal S_1\) already contain the largest weight \(x_1\), yet their sums are still negative. Replacing \(x_1\) by a uniformly chosen unselected weight can only make these sums more negative. On the other hand, since the total weight is nonnegative, the average of all \(k\)-sums must also be nonnegative. Thus, these negative sums must be balanced by positive ones. This suggests that we should consider not only how many negative sums there are, but also their magnitudes.

To measure the magnitudes of negative \(k\)-sums, let \(F\) denote the sum of a uniformly random \(k\)-set, and put \(Y=(-F)_+\), where \(x_+=\max\{x,0\}\). Thus, \(\mathbb E Y\) measures the average magnitude of the negative sums over all \(k\)-sets. When \(x_1\) is exchanged out of a negative set in the star, even a slightly negative sum makes a definite contribution to the negative first moment. We may normalize the weights so that \(\sum_{i=1}^n x_i=0\), \(f(S)\neq0\) for every \(S\in\binom{[n]}{k}\), and \(\mathbb E F_+=\mathbb E Y\) has a fixed positive value, while preserving \(b>0\) and \(c\le b\). The abundance of such negative sets then leaves less room for exceptionally negative \(k\)-sums. Comparing the negative first moments before and after the exchange makes this trade-off precise: with \(s=kx_1/(k-1)\), we obtain an upper bound on \(\mathbb E[(Y-s)_+]\), which measures the average amount by which \(k\)-sums fall below \(-s\).

On the other hand, the positive \(k\)-sums must balance the negative ones, but under our assumption, at most a fraction \(k/n\) of all \(k\)-sets have positive sum. Since normalization fixes \(\mathbb E F_+=\mathbb E Y\), the Cauchy--Schwarz inequality gives
\[
(\mathbb E F_+)^2
\le \Pr(F>0)\,\mathbb E F_+^2
\le \frac{k}{n}\,\mathbb E F_+^2.
\]
Thus, these relatively few positive sums must carry a substantial second moment.

We further sharpen the Cauchy--Schwarz lower bound on \(\mathbb E F_+^2\) by exchanging \(x_1\) into positive sets outside the star. Using a similar exchange argument, we also derive an upper bound on \(\mathbb E F_+^2\) in terms of \(x_1\) and the total second moment \(\mathbb E F^2\). Comparing these two bounds on \(\mathbb E F_+^2\), together with the identity
\[
\mathbb E F^2=\mathbb E F_+^2+\mathbb E Y^2,
\]
yields a lower bound on \(\mathbb E Y^2\).
 This suggests that some negative \(k\)-sums must be quite large in magnitude, seemingly at odds with our earlier bound on the average amount by which \(k\)-sums fall below \(-s\).

However, these two estimates do not directly contradict each other, since a few exceptionally negative sums could make a large contribution to the second moment. To control this possibility, we use Hölder's inequality:
\[
\mathbb E Y^2-s\mathbb E Y
\le(\mathbb E Y^4)^{1/3}
\bigl(\mathbb E[(Y-s)_+]\bigr)^{2/3}.
\]
We therefore seek an upper bound on \(\mathbb E Y^4\). To obtain it, we compute the total fourth moment \(\mathbb E F^4\) and derive a lower bound on \(\mathbb E F_+^4\) from the positive second moment. Subtracting the latter from the former gives an upper bound on \(\mathbb E Y^4\). Combining this with Hölder's inequality and the preceding estimates yields a contradiction, proving Theorem~\ref{thm:main}. To establish Theorem~\ref{thm:stability}, we further show, using the same counting and moment estimates, that \(c/b\) is bounded away from \(1\) by an absolute constant whenever \(b>0\).
\section{Counting negative \texorpdfstring{$k$}{k}-sets via spectral estimates}
\label{sec:complement}

Recall that $b$ counts the negative $k$-sets in $\mathcal S_1$,
and $c$ the nonnegative $k$-sets outside $\mathcal S_1$.
Since $|\mathcal G|=|\mathcal S_1|-b+c$, the condition
$|\mathcal G|\le|\mathcal S_1|$ is equivalent to $c\le b$.
We show that, when $k\ge7$ and $4k<n<5k$, this condition together
with $b>0$ forces a substantial proportion of the star
to be negative. More precisely, we derive lower bounds
on $\beta=b/|\mathcal S_1|$.

For $N\ge k$, let $g(N,k)$ be the minimum number of nonnegative $k$-sums
among all lists of $N$ real weights with nonnegative sum. Let $\mathcal B$
consist of the $(k-1)$-sets obtained by deleting $1$ from the negative
members of $\mathcal S_1$, so $|\mathcal B|=b$. For each $B\in\mathcal B$,
we have $f(B\cup\{1\})<0$. Since the total weight is nonnegative, the
complement of $B\cup\{1\}$ has positive total weight. If $n\ge2k$, this
complement contains at least $g(n-k,k)$ nonnegative $k$-sets, all outside
$\mathcal S_1$ and disjoint from $B$. Conversely, each of the $c$ nonnegative
$k$-sets outside $\mathcal S_1$ is disjoint from at most
$\binom{n-k-1}{k-1}$ members of $\mathcal B$. Counting these pairs gives
\begin{equation}\label{eq:joint}
 b\,g(n-k,k)\le c\binom{n-k-1}{k-1}.
\end{equation}

The inequality \eqref{eq:joint} alone does not give a
lower bound on $b/|\mathcal S_1|$ when $c\le b$.
The maximality of $x_1$ provides an additional
constraint. For every $B\in\mathcal B$ and $j\notin B$,
\[
f(B\cup\{j\})\le f(B)+x_1<0.
\]
Thus, for each $B\in\mathcal B$, at least $g(n-k,k)$
members of $\mathcal G\setminus\mathcal S_1$ are disjoint
from $B$, whereas none contains $B$. We exploit these two properties
using averaging operators.

Let $N\ge2k$. For a function $\varphi$ on the
$k$-subsets of $[N]$, let $P\varphi$ be its average
over the $k$-sets disjoint from its argument.
For a $(k-1)$-set $R$, let $(\mathcal K\varphi)(R)$
and $(\mathcal U\varphi)(R)$ denote the averages of
$\varphi$ over the $k$-sets disjoint from $R$ and
containing $R$, respectively.

In the application on $[n]\setminus\{1\}$, the preceding
observation gives $(\mathcal U\varphi)(B)=0$ for every
$B\in\mathcal B$ whenever $\varphi$ is supported on
$\mathcal G\setminus\mathcal S_1$. Thus, when summing
$(\mathcal K\varphi)(B)$ over $B\in\mathcal B$, we may replace
$\mathcal K$ by $\mathcal K+\omega\mathcal U$ without changing
the sum.

To estimate this modified operator, we use
$\mathcal K=\mathcal UP$. Indeed, choosing a uniform
one-point extension of a $(k-1)$-set $R$, followed
by a uniform $k$-set disjoint from the extension,
produces a uniform $k$-set disjoint from $R$,
since each such set arises from exactly
$N-2k+1$ extensions.

Equip the $k$-set and $(k-1)$-set layers with
their uniform inner products $\langle u,v\rangle=\E[uv]$,
where each expectation is taken uniformly over the corresponding layer.
Write $\|\cdot\|_0$ for the operator norm on mean-zero functions.
Each operator fixes constants and preserves uniform means: averaging over
all its arguments counts every $k$-set equally often. Thus
\[
\mathcal K+\omega\mathcal U=\mathcal U(P+\omega I),\qquad
\|\mathcal K+\omega\mathcal U\|_0
\le\|\mathcal U\|_0\|P+\omega I\|_0.
\]
We estimate the two factors using the eigenvalues
of $\mathcal U^*\mathcal U$ and $P$, respectively.

For $H\subseteq[N]$, define
$e_H(T)=\one_{H\subseteq T}$ on either layer. Fix $|H|=j$, with $1\le j\le k$.
For a $k$-set $T$, $(Pe_H)(T)$ is the probability
that a uniformly chosen $k$-set disjoint from $T$
contains $H$. This probability is zero if
$H\cap T\ne\varnothing$, and otherwise equals
$q_j=(k)_j/(N-k)_j$, where
$(s)_j=s(s-1)\cdots(s-j+1)$ and $(s)_0=1$.
Thus, by inclusion--exclusion,
\[
(Pe_H)(T)
=q_j\one_{H\cap T=\varnothing}
=q_j\sum_{J\subseteq H}(-1)^{|J|}e_J(T).
\]

For $\mathcal U$, put $m=N-k+1$, the number of
one-point extensions of a $(k-1)$-set $R$.
If $R$ contains $H$, every extension contains $H$;
if $R$ misses exactly one element of $H$, exactly
one extension contains $H$; otherwise none does.
Consequently,
\[
\mathcal U e_H
=(1-j/m)e_H+
\frac1m\sum_{h\in H}e_{H\setminus\{h\}}.
\]
The adjoint $\mathcal U^*$ averages over one-point
deletions. Of the $k$ deletions of a $k$-set
containing $H$, exactly $k-j$ still contain $H$.
Hence $\mathcal U^*e_H=(1-j/k)e_H$.

Both calculations express the action of the
operators in terms of $e_H$ and indicators
indexed by proper subsets of $H$.
To separate these lower-order terms, let $V_j$
be the span of the functions $e_H$ with $|H|\le j$
on the $k$-set layer, and put
$W_j=V_j\cap V_{j-1}^{\perp}$.

The preceding identities give
\[
Pe_H\equiv(-1)^jq_je_H,\qquad
\mathcal U^*\mathcal Ue_H\equiv(1-j/k)(1-j/m)e_H
\pmod{V_{j-1}}.
\]
Since $P$ and $\mathcal U^*\mathcal U$ are
self-adjoint and preserve each $V_j$, they also
preserve $W_j$. Projecting the identities above
onto $W_j$ shows that both operators act as
scalar multiples of the identity on $W_j$.
Their eigenvalues there are therefore
\begin{equation}\label{eq:disjoint-levels}
(-1)^jq_j,\qquad (1-j/k)(1-j/m)
\qquad (1\le j\le k).
\end{equation}

Finally, $V_0$ consists of constants and $V_k$
is the full space of functions on $k$-sets.
Thus the spaces $W_j$ span the mean-zero subspace orthogonally.
Since $\|\mathcal U u\|^2=\langle u,\mathcal U^*\mathcal Uu\rangle$,
we obtain $\|\mathcal U\|_0=\sqrt{(1-1/k)(1-1/m)}$.
The eigenvalues of $P$ similarly give
$\|P+\omega I\|_0=\max_{1\le j\le k}|(-1)^jq_j+\omega|$.
On $W_1$, the eigenvalue of $P$ is $-q_1$.
For $k\ge2$, its eigenvalues on $W_j$, $j\ge2$, have
absolute value at most $q_2$, since $N\ge2k$ makes the
$q_j$ nonincreasing. Hence, for $0\le\omega\le q_1$,
\[
\|P+\omega I\|_0
\le\max\{q_1-\omega,\ q_2+\omega\}.
\]
Increasing $\omega$ decreases the first bound
but increases the second. In the proposition below, we choose
$\omega$ so that the second does not exceed the first.

These spectral estimates give a bound relating $g(n-k,k)$
to $\beta$, which also applies when $c>b$, as needed in the proof
of Theorem~\ref{thm:stability}.

\begin{proposition}\label{prop:pair-count}
Suppose $n\ge4k$, $k\ge2$, and $b>0$.
Put $r=k/(n-k)$ and
$\beta=b/\binom{n-1}{k-1}$. Then
\begin{equation}\label{eq:count-linear}
\frac{b}{\max\{b,c\}}\,
\frac{g(n-k,k)}{\binom{n-k-1}{k-1}}
\le \frac{1+r}{8}
\bigl[1+4r+(8-9r)\beta\bigr].
\end{equation}
\end{proposition}

\begin{proof}
All operators now act on $[n]\setminus\{1\}$, so $N=n-1$ and
$m=n-k\ge3k$. Here $q_1=k/(m-1)\ge r$ and
$q_2=k(k-1)/[(m-1)(m-2)]\le r^2$; the latter inequality is
equivalent to $m(m-3k)+2k\ge0$.
Choose $\omega=r(1-r)/2$, so that
$q_2+\omega\le r^2+\omega=r-\omega\le q_1-\omega$.
The spectral estimate above gives $\|P+\omega I\|_0\le q_1-\omega$.
Also, the arithmetic--geometric mean inequality and $m\ge3k$ give
$\|\mathcal U\|_0\le1-1/(2k)-1/(2m)\le1-2/m$.
Using $(1-1/m)q_1=r$ and $2\omega\le q_1$, we obtain
\[
\|\mathcal K+\omega\mathcal U\|_0
\le(1-2/m)(q_1-\omega)\le r-\omega=\frac{r(1+r)}2.
\]

Let $h$ be the indicator of the nonnegative $k$-sets in
$[n]\setminus\{1\}$, with mean $\gamma=c/\binom{n-1}{k}$.
For $B\in\mathcal B$, the operator $\mathcal K$ averages $h$ over
all $\binom{m}{k}$ $k$-subsets of $[n]\setminus(B\cup\{1\})$.
Since $h$ is $1$ on the nonnegative sets and $0$ on the others,
this average is the proportion of these sets that are nonnegative:
\[
(\mathcal K h)(B)
=\frac{|\{G\in\mathcal G\setminus\mathcal S_1:B\cap G=\varnothing\}|}
       {\binom{m}{k}}
\ge\frac{g(m,k)}{\binom{m}{k}}.
\]
The family $\mathcal B$ occupies a proportion $\beta$ of the
$(k-1)$-sets. Averaging over this layer and using
$(\mathcal U h)(B)=0$ for $B\in\mathcal B$ gives
\[
\frac{\beta g(m,k)}{\binom{m}{k}}
\le\langle\one_{\mathcal B},\mathcal K h\rangle
=\langle\one_{\mathcal B},(\mathcal K+\omega\mathcal U)h\rangle.
\]

To apply the mean-zero norm estimate, we subtract the means $\beta$
and $\gamma$ from $\one_{\mathcal B}$ and $h$, respectively.
The operators $\mathcal K$ and $\mathcal U$ fix constants and preserve
uniform means, so $(\mathcal K+\omega\mathcal U)h$ is the sum of
the constant $(1+\omega)\gamma$ and the mean-zero function
$(\mathcal K+\omega\mathcal U)(h-\gamma)$.
Expanding against $\one_{\mathcal B}=\beta+(\one_{\mathcal B}-\beta)$,
the two mixed terms vanish because constants are orthogonal to
mean-zero functions. Hence
\[
\langle\one_{\mathcal B},(\mathcal K+\omega\mathcal U)h\rangle
=(1+\omega)\beta\gamma+
\langle\one_{\mathcal B}-\beta,(\mathcal K+\omega\mathcal U)(h-\gamma)\rangle.
\]
The squared norms of $\one_{\mathcal B}-\beta$ and $h-\gamma$ are
$\beta(1-\beta)$ and $\gamma(1-\gamma)$, respectively.
Cauchy--Schwarz and the operator norm bound therefore yield
\[
\frac{\beta g(m,k)}{\binom{m}{k}}
\le(1+\omega)\beta\gamma+
\frac{r(1+r)}2\sqrt{\beta(1-\beta)\gamma(1-\gamma)}.
\]

Since $\gamma=(c/b)r\beta$, we distinguish two cases.
If $c\le b$, then $\gamma\le r\beta\le1/3$.
Dividing by $r\beta$ and using
$\gamma(1-\gamma)\le r\beta(1-r\beta)$ gives the bound below
with prefactor $1$. If $c>b$, divide instead by $\gamma$ and use
$(1-\gamma)/\gamma\le(1-r\beta)/(r\beta)$; the prefactor is then
$r\beta/\gamma=b/c$.
Thus, using $\binom{m-1}{k-1}=r\binom{m}{k}$, both cases give
\[
\frac{b}{\max\{b,c\}}\,
\frac{g(m,k)}{\binom{m-1}{k-1}}
\le(1+\omega)\beta+
\frac{1+r}{2}\sqrt{r(1-\beta)(1-r\beta)}.
\]
Finally, apply $2\sqrt{xy}\le2x+y/2$ with
$x=r(1-\beta)$ and $y=1-r\beta$.
Substituting $\omega=r(1-r)/2$ bounds the right-hand side by
$(1+r)[1+4r+(8-9r)\beta]/8$, proving \eqref{eq:count-linear}.
\end{proof}

For $0\le\beta\le1$, the right-hand side of
\eqref{eq:count-linear} is increasing in both $r\in[0,1/3]$
and $\beta$. Since $r\le1/3$, we have
\begin{equation}\label{eq:count-uniform}
\frac{b}{\max\{b,c\}}\,
\frac{g(n-k,k)}{\binom{n-k-1}{k-1}}
\le\frac7{18}+\frac56\beta.
\end{equation}

To obtain a lower bound on $\beta$ from
\eqref{eq:count-linear}, we next bound $g(n-k,k)$.

\begin{lemma}\label{lem:weak-boundary}
For $k\ge2$, we have $g(3k+1,k)\ge\tfrac78\binom{3k-2}{k}$.
\end{lemma}
\begin{proof}
Motivated by the example with three weights equal to
$-(3k-2)$ and the remaining weights equal to $3$,
whose nonnegative $k$-sets are precisely those avoiding
the three negative entries, we mark the indices $3k-1$, $3k$ and $3k+1$, which correspond
to the three smallest weights. If every $k$-subset of $[3k-2]$
is nonnegative, then $|\mathcal G|\ge\binom{3k-2}{k}$, as required.
Otherwise, choose a $k$-subset $S$ of $[3k-2]$ with $f(S)<0$.
By the ordering of the weights, each marked index $i$ satisfies
\[
x_i\le \min_{s\in S}x_s\le \frac{f(S)}{k}<0.
\]
Thus deleting any marked index removes a negative weight from a
list of nonnegative total weight. The remaining $3k$ weights
therefore have positive sum, so every partition of the remaining
indices into three $k$-sets has at least one block of positive
weight.

For $j=0,1,2$, let $p_j$ be the proportion of nonnegative sets
among the $k$-sets containing exactly $j$ of the three marked
indices. After one marked index is deleted, the other two may lie
in different blocks of the partition or in the same block. We
consider these two cases separately. In each case, choose the
index to be deleted uniformly from the three marked indices,
then choose uniformly among the ordered partitions of the remaining
$3k$ indices into three $k$-sets satisfying the specified condition.

In either case, the random deletion and the choice of partition
treat the three marked indices symmetrically, and likewise the
unmarked indices. Thus, conditional on containing exactly $j$
marked indices, a given block is uniform among all such $k$-sets
and is nonnegative with probability $p_j$.

When the two remaining marked indices lie in different blocks,
one block is unmarked and the other two have one mark each, so
the expected number of nonnegative blocks is $p_0+2p_1$.
When they lie in the same block, the other two blocks are
unmarked, so the corresponding expectation is $2p_0+p_2$.
Since every partition has at least one nonnegative block, we obtain
\[
p_0+2p_1\ge1,\qquad 2p_0+p_2\ge1.
\]

There are $\binom{3k-2}{k}$ unmarked $k$-sets. Relative to this
number, the numbers of $k$-sets with exactly one and exactly two
marked indices are, respectively,
\[
\frac{3\binom{3k-2}{k-1}}{\binom{3k-2}{k}}
=\frac{3k}{2k-1}\ge\frac32,\qquad
\frac{3\binom{3k-2}{k-2}}{\binom{3k-2}{k}}
=\frac{3(k-1)}{2(2k-1)}\ge\frac12,
\]
where the last inequality uses $k\ge2$. Counting only the
nonnegative $k$-sets with at most two marked indices, we obtain
\[
\frac{|\mathcal G|}{\binom{3k-2}{k}}
\ge p_0+\frac32p_1+\frac12p_2
\ge \frac34(p_0+2p_1)+\frac18(2p_0+p_2)
\ge \frac78.\qedhere
\]
\end{proof}
Dividing the lower bound in Lemma~\ref{lem:weak-boundary} by
$\binom{3k}{k-1}$, we obtain, for $k\ge2$,
\begin{equation}\label{eq:boundary-ratio}
 \frac{g(3k+1,k)}{\binom{3k}{k-1}}
 \ge \frac{7(8k^2-2)}{24k(3k-1)}
 > \frac{41}{54}+\frac{5}{24k}.
\end{equation}
For the strict inequality, subtracting the rightmost expression
from the middle one gives $(12k^2+29k-81)/[216k(3k-1)]$,
which is positive for $k\ge2$.

\begin{proposition}\label{prop:three-quarters}
If $3k\le N\le4k$, then
\(
 g(N,k)>\frac{41}{54}\binom{N-1}{k-1}.
\)
\end{proposition}

For $k=1$, at least one weight is nonnegative, so $g(N,1)\ge1$.
Assume now that $k\ge2$. At $N=3k$ and $N=4k$, averaging over
partitions into $k$-sets gives the stronger bound
$g(N,k)\ge\binom{N-1}{k-1}$. The case $N=3k+1$ follows from
\eqref{eq:boundary-ratio}. In the range $3k+2\le N<4k$, the full
star bound holds for $k=3$~\cite{Chow} and for $4\le k\le7$~\cite{HS}.
This range is empty when $k=2$.

It remains to consider $k\ge8$ and $3k+2\le N<4k$.
We defer this case to Appendix~\ref{app:moment-comparison},
where \eqref{eq:joint} is combined with the moment estimates
developed in Section~\ref{sec:tools}. These estimates are
established independently of both Proposition~\ref{prop:three-quarters}
and Lemma~\ref{lem:countgate}.

\begin{lemma}\label{lem:countgate}
Suppose $k\ge7$ and $4k<n<5k$. If $\mathcal S_1$ has a negative member
and $|\mathcal G|\le\binom{n-1}{k-1}$, then the fraction $\beta$ of
negative sets in $\mathcal S_1$ satisfies
\[
 \frac{n-1}{k}\beta>\frac{16}9+\frac1k .
\]
\end{lemma}
\begin{proof}
Suppose first that $n=4k+1$. Since $b/\max\{b,c\}=1$ under the
hypotheses of the lemma, \eqref{eq:boundary-ratio} and
\eqref{eq:count-uniform} give
\[
 \frac{41}{54}+\frac{5}{24k}
 < \frac{g(3k+1,k)}{\binom{3k}{k-1}}
 \le \frac7{18}+\frac56\beta.
\]
Thus $\beta>4/9+1/(4k)$, which gives the desired bound since
$(n-1)/k=4$.

Now suppose that $n\ge4k+2$, so that $r\le k/(3k+2)$. If the
desired bound were false, then
\[
 \beta\le\frac{16k+9}{9(n-1)}
 \le\frac{16k+9}{9(4k+1)}.
\]
Since the right-hand side of \eqref{eq:count-linear} is increasing
in both $r$ and $\beta$, these bounds give
\[
 \frac{g(n-k,k)}{\binom{n-k-1}{k-1}}
 \le \frac{41}{54}
 -\frac{21k^2-127k-79}{54(4k+1)(3k+2)^2}
 < \frac{41}{54},
\]
where the strict inequality follows from
$21k^2-127k-79=(k-7)(21k+20)+61>0$ for $k\ge7$.
This contradicts Proposition~\ref{prop:three-quarters} applied
with $N=n-k$, since $3k<n-k<4k$.
\end{proof}

\section{Moment inequalities for \texorpdfstring{$k$}{k}-sums}
\label{sec:tools}

For $k\ge7$ and $4k<n<5k$, Lemma~\ref{lem:countgate}
gives a lower bound on the proportion of negative sets
in the star under the assumptions $b>0$ and $c\le b$.
We now use moment estimates to show that these
assumptions are incompatible.

We derive the moment estimates for a list of $N$ weights, since we shall
also need them on the complement of a negative $k$-set. Let
$x_1\ge\cdots\ge x_N$ have nonnegative sum, and let
$F=\sum_{i\in S}x_i$ for a uniformly chosen $k$-subset $S$ of $[N]$,
where $1\le k<N$. On a complement we relabel the remaining indices in
their original order, so that $\mathcal S_1$ is always the star at an
index of largest weight.

We first reduce to the case $x_1+\cdots+x_N=0$ with $f(S)\ne0$ for
every $k$-set $S$. Subtracting the mean $\bar x=N^{-1}\sum_i x_i$
from each weight changes $f(S)$ to $f(S)-k\bar x\le f(S)$.
Now add $(N-1)\delta$ to $x_1-\bar x$ and subtract $\delta$ from
each $x_j-\bar x$ with $j>1$. The total remains zero.
A $k$-sum increases by $(N-k)\delta$ if its set contains $1$ and
decreases by $k\delta$ otherwise. Since there are finitely many
$k$-sets, we can take $\delta>0$ small enough to preserve the signs
of all nonzero centered sums. Every zero sum becomes nonzero, $x_1$
is now the unique largest weight, and every originally negative set
remains negative. Thus $b$ can only increase and $c$ can only decrease.

Retain the notation $x_j$ and $F$ for the adjusted weights and sums.
Put $p=k/N$, $q=1-p$ and $d=\sum_jx_j^2$. In
$F^2=\sum_{j\in S}x_jF$, the term $x_jF$ occurs precisely when
$j\in S$. Taking expectations gives
\begin{equation}\label{eq:second-moment-expansion}
 \E F^2=\sum_{j=1}^N x_j\,\E[F\one_{j\in S}].
\end{equation}
Given $j\in S$, each $i\ne j$ is selected with probability
$(k-1)/(N-1)$. Since $\sum_{i\ne j}x_i=-x_j$, we obtain
\[
 \E[F\one_{j\in S}]
 =p\,\E[F\mid j\in S]
 =p\left(x_j-\frac{k-1}{N-1}x_j\right)
 =\lambda x_j,\qquad \lambda=\frac{k(N-k)}{N(N-1)}.
\]
Substitution in \eqref{eq:second-moment-expansion} gives
$\E F^2=\lambda d$.

Write $x_+=\max\{x,0\}$, $x_-=\max\{-x,0\}$ and $Y=F_-$.
Since $\E F=p\sum_jx_j=0$ and no $k$-sum vanishes, $\E F_+=\E Y>0$.
Rescale all weights by the same positive factor so that
$\E F_+=\E Y=\lambda$. We call these weights \emph{normalized}.
The rescaling preserves every sign, and $|x_j|\le1$ for all $j$,
because
\[
 -\lambda=-\E Y\le\E[F\one_{j\in S}]=\lambda x_j\le\E F_+=\lambda .
\]
Also $x_1>0$, since the weights sum to zero and are not all zero.

\subsection*{A first-moment bound}

Choose $S_s$ uniformly from $\mathcal S_1$, then choose $j$ uniformly
from $[N]\setminus S_s$. Replacing $1$ by $j$ gives a set
$S_o=(S_s\setminus\{1\})\cup\{j\}$ outside the star. Every $k$-set
$U$ avoiding $1$ arises from exactly $k$ equally likely pairs,
one for each $j\in U$, with $S_s=(U\setminus\{j\})\cup\{1\}$.
Thus $S_o$ is uniform outside the star and, conditional on $S_o$,
$j$ is uniform in $S_o$.

Write $F_s=f(S_s)$ and $F_o=f(S_o)$. We have $F_s-F_o=x_1-x_j\ge0$.
Given $S_o$, the mean of $x_j$ is $F_o/k$, so
\begin{equation}\label{eq:coupled-mean}
 \E[F_s\mid S_o]=x_1+(1-v)F_o,\qquad v=1/k .
\end{equation}
Let $\beta=\Pr(F_s<0)$ and $\widehat\beta=(N-1)\beta/k$. For $k\ge2$,
the conditional mean in \eqref{eq:coupled-mean} is negative precisely
when $F_o<-kx_1/(k-1)$.

\begin{lemma}\label{lem:first-moment}
For normalized weights,
\begin{equation}\label{eq:star-slack}
 x_1\widehat\beta\le1 .
\end{equation}
If $k\ge2$ and $s=kx_1/(k-1)$, then
\begin{equation}\label{eq:slack}
 (1-v)\E[(Y-s)_+]
 \le\frac{N-1}{N}\E[(F_s)_-]\le\lambda(1-x_1\widehat\beta).
\end{equation}
\end{lemma}
\begin{proof}
Replacing $1$ by $j$ can only increase the negative part, and on
$F_s<0$ the increase is exactly $x_1-x_j$. Given $S_s$, the mean
of $x_j$ is $-F_s/(N-k)$, since the weights outside $S_s$ sum to
$-F_s$. Consequently
\[
 \E[(F_o)_-]-\E[(F_s)_-]
 \ge\E[(x_1-x_j)\one_{F_s<0}]
 =x_1\beta-\frac{\E[(F_s)_-]}{N-k}.
\]
Multiplying by $q$ and using
$p\E[(F_s)_-]+q\E[(F_o)_-]=\lambda$ and
$q\beta=\lambda\widehat\beta$ gives
\[
 \frac{N-1}{N}\E[(F_s)_-]\le\lambda(1-x_1\widehat\beta).
\]
This also proves \eqref{eq:star-slack}.

For $k\ge2$, taking the negative part in \eqref{eq:coupled-mean}
and applying conditional Jensen's inequality gives
\[
 (1-v)\E[((F_o)_--s)_+]
 =\E\big[(\E[F_s\mid S_o])_-\big]\le\E[(F_s)_-].
\]
Inside the star we have $((F_s)_--s)_+\le(F_s)_-$. Multiply the
outside estimate by $q$ and the inside estimate by $(1-v)p$.
Adding them gives \eqref{eq:slack}, since $q+(1-v)p=(N-1)/N$.
\end{proof}

\subsection*{The second moment}

For normalized weights put $d_+=\E F_+^2/\lambda$,
$d_-=\E Y^2/\lambda$ and $\mu=\Pr(F>0)$, so that $d=d_++d_-$.
The conditional mean of $F$ given $F>0$ is $\lambda/\mu$, and hence
\begin{equation}\label{eq:positive-second-variance}
 d_+=\frac{\lambda}{\mu}
 +\frac1\lambda\E\big[(F-\lambda/\mu)^2\one_{F>0}\big].
\end{equation}
Exchanging $1$ into positive sets outside the star gives a lower
bound on $\E[(F-h)^2\one_{F>0}]$ for any $h\in\mathbb R$. At
$h=\lambda/\mu$ this refines the Cauchy--Schwarz bound
$d_+\ge\lambda/\mu$. For the upper bound on $d_+$ we use exchanges
between sets of opposite signs. Write $r=k/(N-k)$ and $\gamma=\Pr(F_o>0)$.
Since no $k$-sum vanishes,
$\mu=p(1-\beta)+q\gamma$.

\begin{lemma}\label{lem:exchange-bound}
For normalized weights and any $h\in\mathbb R$,
\begin{equation}\label{eq:conditional-x}
 d_+\ge2h-\mu h^2/\lambda+\frac{N-1}{N}\gamma(x_1-vh)^2 .
\end{equation}
For every $H\ge\max\{x_1,(k-1)/(N-k)\}$,
\begin{equation}\label{eq:refex}
 d_+\le\frac{N-k}{N-1}H+
        \left(\mu+\frac{r}{(N-1)H}\right)d .
\end{equation}
\end{lemma}
\begin{proof}
If $F_o>0$, then $F_s>0$ as well. To compare the deviations of these
two sums from $h$, write \eqref{eq:coupled-mean} as
\[
 \E[F_s-h\mid S_o]-(1-v)(F_o-h)=x_1-vh .
\]
Weighted Cauchy--Schwarz, followed by conditional Jensen, gives
\[
 pq(x_1-vh)^2\le p\E[(F_s-h)^2\mid S_o]+q(F_o-h)^2,
\]
since $q+p(1-v)^2\le1$. Integrate over $F_o>0$. In the first
term we may enlarge the event to $F_s>0$, so the two terms
together give
\[
 pq\gamma(x_1-vh)^2\le\E[(F-h)^2\one_{F>0}]
 =\lambda d_+-2h\lambda+\mu h^2 .
\]
As $pq/\lambda=(N-1)/N$, this proves \eqref{eq:conditional-x}.

For \eqref{eq:refex}, observe that
$\lambda(d_+-\mu d)=\E[(F^2-\lambda d)\one_{F>0}]$.
We express $F^2-\lambda d$ as a sum over exchanges out of $S$.
For $j\in S$ and $l\notin S$, let
$S'=(S\setminus\{j\})\cup\{l\}$, $t=f(S\setminus\{j\})$ and
$\Delta=x_j-x_l=f(S)-f(S')$. Multiply $\Delta$ by $t+q(x_j+x_l)$.
This factor is unchanged when the exchange is reversed, and the
choice of $q$ gives the two sums of squares below the same coefficient.

Fix $S$ and write $F=f(S)$. Summing over $j\in S$ and $l\notin S$
gives $\sum_{j,l}\Delta=NF$ and $\sum_{j,l}x_jx_l=-F^2$. Since
$t+q(x_j+x_l)=F-px_j+qx_l$, expansion gives
\[
 \begin{aligned}
 \sum_{j,l}\Delta[t+q(x_j+x_l)]
 &=(N-1)F^2-p(N-k)\sum_{j\in S}x_j^2-qk\sum_{l\notin S}x_l^2\\
 &=(N-1)(F^2-\lambda d),
 \end{aligned}
\]
where $p(N-k)=qk=(N-1)\lambda$.

Now sum the identities for $\Delta$ and $\Delta[t+q(x_j+x_l)]$
over all positive sets $S$.
An exchange between two positive sets occurs in both directions.
Reversing it changes the sign of $\Delta$ and preserves
$t+q(x_j+x_l)$, so these terms cancel in both sums. Over the
remaining exchanges, all of which satisfy $f(S)>0>f(S')$, we have
$\Delta>0$ and $\sum\Delta=N\binom Nk\E F_+=N\binom Nk\lambda$.
Dividing the weighted sum by this total gives
\begin{equation}\label{eq:exchange-id}
 \frac{N-1}{N}(d_+-\mu d)=
 \frac{\sum\Delta[t+q(x_j+x_l)]}{\sum\Delta}.
\end{equation}

Group the exchanges in \eqref{eq:exchange-id} by
$R=S\cap S'$ and the inserted index $l$. In one such group,
$t=f(R)$ and $y_i=t+x_i$ for $i\notin R$. Thus $y_l<0$ and the
possible starting sets are $R\cup\{j\}$ with $y_j>0$. Let $\phi$
be the average of $t+q(x_j+x_l)$ in this group with weights
$\Delta=y_j-y_l$. The ratio in \eqref{eq:exchange-id} is a weighted
average of these group means, so \eqref{eq:refex} follows if each
satisfies $\phi\le qH+pd/[(N-k)H]$. Since $x_j\le H$,
\[
 t+q(x_j+x_l)=pt+q(x_j+y_l)\le pt+qH .
\]
Hence $\phi-qH\le pt$. We need only consider $\phi>qH$, so $t>0$.

Using $t+q(x_j+x_l)=(1-2q)t+q(y_j+y_l)$, we obtain
\[
 \phi=(1-2q)t+
 q\,\frac{\sum_{y_j>0}(y_j^2-y_l^2)}
              {\sum_{y_j>0}(y_j-y_l)} .
\]
The denominator is at least
$\sum_{y_j>0}y_j\ge\sum_{j\notin R}y_j=(N-k)t$.
The numerator is at most $\sum_{j\notin R}y_j^2$. Since
$\sum_{j\notin R}x_j=-t$,
\[
 \sum_{j\notin R}y_j^2=(N-k-1)t^2+\sum_{j\notin R}x_j^2
 \le (N-k-1)t^2+d .
\]
Using $Nq=N-k$ gives $Nt\phi\le(k-1)t^2+d$. Since $y_l<0$ and
$|x_l|\le1$, we have $t<-x_l\le1$. Also $(N-k)H\ge k-1$.
Thus the bracket in the first line below
is nonnegative, and we may use $pt\ge\phi-qH$ to obtain
\[
 \begin{aligned}
 pd&\ge pNt(\phi-qH)+pt\big[(N-k)H-(k-1)t\big]\ge(\phi-qH)\big[(N-k)H+t\big].
 \end{aligned}
\]
This proves the bound on $\phi$. Averaging over the groups in
\eqref{eq:exchange-id} proves \eqref{eq:refex}.
\end{proof}

Write $P=\mu+r/[(N-1)H]$ for the coefficient of $d$ in \eqref{eq:refex}. Substituting
$d=d_++d_-$ gives $P d_-\ge(1-P)d_+-(N-k)H/(N-1)$. If $P<1$,
we may insert \eqref{eq:conditional-x} to obtain
\begin{equation}\label{eq:negative-second}
 P d_-\ge(1-P)\left(2h-\frac{\mu h^2}{\lambda}
       +\frac{N-1}{N}\gamma(x_1-vh)^2\right)-\frac{N-k}{N-1}H .
\end{equation}
For the comparison with H\"older's inequality we shall retain both bounds
on $d_+$: together they bound $d$ from below, while the upper bound alone
gives $d_-=d-d_+$ a lower bound that still depends on $d$.

We also consider $w=d_+/d$, the fraction of the total second moment
contributed by positive sums. Since $d_-=d_+(1/w-1)$, an upper
bound on $w$ can be combined with a lower bound on $d_+$ to bound
$d_-$ from below. If $\mu\le tp$, then $d_+\ge\lambda/\mu$ gives
$1/d\le tpw/\lambda$. Dividing \eqref{eq:refex} by $d$ and using
$(N-k)/(N-1)=\lambda/p$ yields
\begin{equation}\label{eq:w-range}
 w\le\frac{tp+r/[(N-1)H]}{1-tH}
 \qquad\bigl(H\ge\max\{x_1,(k-1)/(N-k)\},\ tH<1\bigr).
\end{equation}

To obtain a lower bound on $w$, we average over disjoint $k$-sets.
In the final moment comparison, this bound and the exchange estimates
give $x_1\ge (k-1)/(N-k)$, so we may take $H=x_1$ in \eqref{eq:refex}.
The same averaging argument also gives a lower bound on the proportion
of nonnegative $k$-sums. In Appendix~\ref{app:moment-comparison},
we apply this second bound to the weights outside a negative $k$-set.
\begin{proposition}\label{prop:half}
For normalized weights with $N>2k$ and $\mu\le tp$, the fraction
$w=d_+/d$ satisfies
\begin{equation}\label{eq:disjoint-w}
 1+(1+r)w^2\le[2+r+t(1-r)]w .
\end{equation}
For any $N$ real weights of nonnegative total and $1\le k<N$, the
proportion of nonnegative $k$-sums is at least
$\min\{k,N-k\}/[2(N-1)]$.
\end{proposition}
\begin{proof}
First assume that the weights are normalized and $N>2k$. For a
$k$-set $S$, let $(PG)(S)$ be the average of $G$ over the $k$-sets
disjoint from $S$. Each index outside $S$ belongs to a fraction $r$
of these sets. Since the weights outside $S$ sum to $-F$, we have
$PF=-rF$. By \eqref{eq:disjoint-levels}, the eigenvalues of $P$
on the mean-zero subspace are at most
$e=k(k-1)/[(N-k)(N-k-1)]<r$.

Take $G=F_+$ or $G=Y$, and put $u=\E G^2/(\lambda d)$, so that
$u=w$ or $u=1-w$, respectively. The constant component of $G$ is
$\lambda$. Its projection onto the span of $F$ is $uF$ in the first
case and $-uF$ in the second, since $FG=G^2$ or $FG=-G^2$.
These two components are orthogonal, and the remaining component
has squared norm $\lambda d\,u(1-u)-\lambda^2$. Since $P$ is
self-adjoint, fixes constants and satisfies $PF=-rF$, we obtain
\[
 0\le\langle G,PG\rangle
 \le\lambda^2-ru^2\lambda d
       +e\big[\lambda d\,u(1-u)-\lambda^2\big].
\]
Here the first inequality follows from $G,PG\ge0$. Rearranging
and using $wd=d_+\ge\lambda/\mu$ gives
\[
 u[(r+e)u-e]\le(1-e)\frac{\lambda}{d}\le(1-e)\mu w .
\]
For $G=Y$ we have $u=1-w$. In the bound on $\langle G,PG\rangle$
we may replace $e$ by $r^2$, since $e\le r^2$ and its coefficient
is a squared norm. Substituting $\mu\le tp$ and
$p=r/(1+r)$ gives $(1-w)[1-(1+r)w]\le t(1-r)w$, which expands
to \eqref{eq:disjoint-w}.

For the counting bound, return to
$u[(r+e)u-e]\le(1-e)\mu w$ and take $u=\max\{w,1-w\}$.
Since $u\ge1/2$, we have $(r+e)u-e\ge(r-e)/2>0$.
As $u\ge w>0$, it follows that
$\mu\ge(r-e)/[2(1-e)]=k/[2(N-1)]$.
Every originally negative $k$-set remains negative after normalization,
so the original weights have at least as many nonnegative $k$-sums
as the normalized weights. If $N=2k$, at least one set in each complementary pair
is nonnegative. If $N<2k$, first center the weights. Each $k$-sum
then equals the sum of the negated weights on its complement.
Applying the bound just proved for sets of size $N-k$ gives the
result, since centering can only decrease the nonnegative count.
\end{proof}

When $N=2k+1$, Proposition~\ref{prop:half} gives a lower bound of $1/4$
on the proportion of nonnegative $k$-sums. The minimum is attained
when two weights are equal to $-(2k-1)$ and the remaining $2k-1$
weights are equal to $2$~\cite{BM,MM}. A $k$-set has nonnegative sum
exactly when it avoids both negative weights, so the minimum proportion
is $\binom{2k-1}{k}/\binom{2k+1}{k}=(k+1)/(4k+2)$. This tends to
$1/4$ as $k\to\infty$, so the factor $1/2$ in the counting bound
cannot be increased.

\subsection*{The fourth moment and H\"older's inequality}

The exchanges give a lower bound on $\E Y^2$, whereas
\eqref{eq:slack} bounds $\E[(Y-s)_+]$ from above. A few very large
values of $Y$ could account for the second moment while contributing
little to the first. To compare the two estimates by H\"older's
inequality, we therefore bound $\E Y^4$. We compute $\E F^4$ and
subtract a lower bound on the contribution from $F_+$.

\begin{lemma}\label{lem:fourth-moment}
For normalized weights let $K=\E F^4/(\lambda d^2)$. If $N\ge4$, then
\begin{equation}\label{eq:fourth-mom}
 K=3\theta+(1-6\theta)\frac{\sum_jx_j^4}{d^2},\qquad
 \theta=\frac{(k-1)(N-k-1)}{(N-2)(N-3)} .
\end{equation}
In particular $K\le\max\{3\theta,1-3\theta\}$. For fixed $k\ge2$,
$\theta$ decreases with $N$ on $N\ge2k$. For all $N>k$,
\begin{equation}\label{eq:split}
 \frac{\E Y^4}{\lambda}\le Kd^2-d_+^3,\qquad d\le4K .
\end{equation}
\end{lemma}
\begin{proof}
Suppose $N\ge4$. Write $A=\sum_i x_i^4$ and let $b_j=(k)_j/(N)_j$
be the probability that $S$ contains $j$ prescribed indices. In the
expansion of $F^4$, a term involving $j$ distinct indices is present
with probability $b_j$. We therefore group the terms of
$(\sum_i x_i)^4$ by their number of distinct indices.

The terms involving a single index sum to $A$. Those involving
two distinct indices sum to
\[
 4\sum_{i\ne j}x_i^3x_j+6\sum_{i<j}x_i^2x_j^2
 =-4A+3(d^2-A)=3d^2-7A .
\]
For three distinct indices, fix the repeated index $i$. Each monomial
$x_i^2x_jx_l$ with $j<l$ occurs $12$ times.
Since $\sum_{j<l,\,j,l\ne i}x_jx_l=x_i^2-d/2$, the sum of these
terms is $12\sum_i x_i^2(x_i^2-d/2)=12A-6d^2$. The remaining terms
sum to $3d^2-6A$, as $(\sum_i x_i)^4=0$. Multiplying the four
sums by $b_1,\ldots,b_4$, respectively, gives
\[
 \E F^4=(b_1-b_2)A+3(b_2-2b_3+b_4)(d^2-2A).
\]
Here $b_1-b_2=\lambda$. The coefficient $b_2-2b_3+b_4$ is the
probability of including two prescribed indices and avoiding two
others, namely $(k)_2(N-k)_2/(N)_4=\lambda\theta$. This proves
\eqref{eq:fourth-mom}. The bound on $K$ follows from $0\le A/d^2\le1$.
For fixed $k\ge2$ and $N\ge2k$, monotonicity follows from
\[
 \theta(N)-\theta(N+1)=
 \frac{(k-1)(N-2k+1)}{(N-3)(N-2)(N-1)}>0 .
\]

For every $N>k$, we use H\"older's inequality in the form
$(\E G^2)^3\le(\E G)^2\E G^4$ for $G\ge0$. Applying it to $F_+$ and $Y$,
both of mean $\lambda$, gives $\E F_+^4/\lambda\ge d_+^3$ and
$\E Y^4/\lambda\ge d_-^3$. Subtracting the first lower bound
from $\E F^4/\lambda=Kd^2$ gives the upper bound on $\E Y^4$.
Adding the two lower bounds gives
$Kd^2\ge d_+^3+d_-^3\ge(d_++d_-)^3/4$, whence $d\le4K$.
\end{proof}

Let $k\ge2$ and recall $s=kx_1/(k-1)$ from
Lemma~\ref{lem:first-moment}. Since $0\le(Y-s)_+\le Y$,
H\"older's inequality with exponents $3$ and $3/2$ gives
\[
 \begin{aligned}
 \lambda(d_--s)=\E[Y(Y-s)]
 &\le\E\big[Y^{4/3}(Y-s)_+^{2/3}\big]\le(\E Y^4)^{1/3}\big(\E[(Y-s)_+]\big)^{2/3}.
 \end{aligned}
\]
Inserting \eqref{eq:slack} yields
\begin{equation}\label{eq:tail-interpolation}
 d_--s\le
 \left(\frac{\E Y^4}{\lambda}\right)^{1/3}
 \left(\frac{1-x_1\widehat\beta}{1-v}\right)^{2/3}.
\end{equation}

When $k\ge7$ and $4k\le N\le5k$, Lemma~\ref{lem:fourth-moment} gives
$7/48\le\theta(5k)\le\theta(N)\le\theta(4k)\le3/16$ and hence
$K\le\max\{3\theta,1-3\theta\}\le9/16$. Thus \eqref{eq:split} gives
$\E Y^4/\lambda\le Kd^2-d_+^3\le9d^2/16$.
Substituting in \eqref{eq:tail-interpolation}, we obtain
\begin{equation}\label{eq:holder-upper}
 d_--s\le
 \left(\frac{3d}{4}\,\frac{1-x_1\widehat\beta}{1-v}\right)^{2/3}.
\end{equation}
This bound holds regardless of the number of nonnegative sets. We now
show that if $b>0$ and $c\le b$, the exchange estimates force $d_--s$
to be strictly larger than the right-hand side.

\begin{proposition}\label{prop:strict-star}
Suppose that $k\ge7$ and $4k<n<5k$. If $b>0$, then $c>b$.
\end{proposition}
\begin{proof}
Suppose that $b>0$ and $c\le b$. Normalize the weights with $N=n$.
This can only increase $b$ and decrease $c$, so
$|\mathcal G|=|\mathcal S_1|-b+c\le|\mathcal S_1|$ still holds.
Put $\kappa=16/9+v$. Lemma~\ref{lem:countgate} and
\eqref{eq:star-slack} give $\widehat\beta>\kappa$ and
$x_1\widehat\beta\le1$. In particular, $0<x_1<1/\kappa$.

To apply the exchange estimates, write $T=n/k$, $t=\mu/p$ and
$Q=\lambda/p=(T-1)/(T-v)$. Since no normalized sum vanishes,
$t=|\mathcal G|/|\mathcal S_1|\le1$. We take $h=Q$ in
\eqref{eq:conditional-x}: this maximizes $2h-\mu h^2/\lambda$ when
$\mu=p$. Substituting $q\gamma=\mu-p(1-\beta)$ there and rewriting
\eqref{eq:refex} gives the following bounds, before we use $t\le1$:
\begin{equation}\label{eq:star-second}
 \begin{aligned}
 d_+&\ge Q+\frac{\widehat\beta(x_1-vQ)^2}{T(T-1)}
           +(1-t)\left(Q-\frac{p}{Q}(x_1-vQ)^2\right),\\
 d_+&\le QH+\left(\frac{t}{T}
           +\frac{v}{H(T-1)(T-v)}\right)d
 \qquad\left(H\ge\max\left\{x_1,\frac{1-v}{T-1}\right\}\right).
 \end{aligned}
\end{equation}
The coefficient of $1-t$ is positive, since $Q\ge3/4$, $p\le1/4$ and
$x_1,vQ\in[0,1]$. We may therefore set $t=1$ in both bounds and replace
$\widehat\beta$ by $\kappa$ in the lower bound. Also, since
$r=1/(T-1)>1/4$, \eqref{eq:disjoint-w} with $t=1$ gives
$1+5w^2/4\le3w$, so $w=d_+/d\ge2/5$.

Now $v\le1/7$ and $4+v\le T<5$. Applying
Lemma~\ref{lem:moment-comparison} with $a=x_1$ to these second-moment
bounds gives
\[
 d_--s>\left(\frac{3d}{4}\,\frac{1-\kappa x_1}{1-v}\right)^{2/3}.
\]
But $0\le1-x_1\widehat\beta\le1-x_1\kappa$, so H\"older's bound
\eqref{eq:holder-upper} gives
\[
 d_--s\le\left(\frac{3d}{4}\,\frac{1-x_1\widehat\beta}{1-v}\right)^{2/3}
 \le\left(\frac{3d}{4}\,\frac{1-x_1\kappa}{1-v}\right)^{2/3},
\]
contradicting the strict lower bound.
\end{proof}

\section{Proofs of Theorems~\ref{thm:main} and~\ref{thm:stability}}
\label{sec:main-moment}

Recall that $b=|\mathcal S_1\setminus\mathcal G|$ and
$c=|\mathcal G\setminus\mathcal S_1|$. Thus
\begin{equation}\label{eq:excess}
 |\mathcal G|-|\mathcal S_1|=c-b,\qquad
 |\mathcal G\mathbin\triangle\mathcal S_1|=c+b.
\end{equation}
When $b=0$, both quantities equal $c$. For $b>0$,
Proposition~\ref{prop:strict-star} gives $c>b$ when $k\ge7$ and
$4k<n<5k$. To complete the proof of Theorem~\ref{thm:main}, we use
the divisible case~\cite{BM,ChowThesis,MS} and a reduction from $n$ points to $n-k$ points~\cite{MM}.

\begin{lemma}[{\cite{BM,ChowThesis,MS}}]\label{lem:divisible-star}
If $k\mid n$, then $|\mathcal G|\ge\binom{n-1}{k-1}$. If in addition
$n>2k$, equality holds if and only if $\mathcal G=\mathcal S_1$.
\end{lemma}

To identify the families attaining equality in the reduction, we use the following theorem.

\begin{theorem}[{\cite{EKR,HM}}]\label{thm:ekr}
Let $k\ge1$ and $n\ge2k$, and let $\mathcal F\subseteq\binom{[n]}{k}$
be intersecting, that is, any two distinct members of $\mathcal F$ have
nonempty intersection. Then $|\mathcal F|\le\binom{n-1}{k-1}$.
If $n>2k$, equality holds if and only if there is an $i\in[n]$ such that
$\mathcal F=\{S\in\binom{[n]}{k}:i\in S\}$.
\end{theorem}

\begin{lemma}\label{lem:disjoint-reduction}
Suppose that $n>2k$ and
$g(n-k,k)\ge\binom{n-k-1}{k-1}$.
Then $|\mathcal G|\ge\binom{n-1}{k-1}$,
with equality if and only if $\mathcal G=\mathcal S_1$.
\end{lemma}

\begin{proof}
The lower bound follows directly from the disjoint-pair
reduction of Manickam and Mikl\'os~\cite{MM}, together with
the hypothesis $g(n-k,k)\ge\binom{n-k-1}{k-1}$.
It remains to characterize the equality case.

Suppose that $|\mathcal G|=\binom{n-1}{k-1}$.
Let $(U,V)$ be a uniformly random ordered pair of disjoint
$k$-sets, and put $\alpha=|\mathcal G|/\binom nk=k/n$.
Retaining the contribution from pairs with both members
in $\mathcal G$ in the disjoint-pair counting argument, we obtain
\[
\begin{aligned}
\frac{k}{n}
&\ge \Pr(U\in\mathcal G,\ V\in\mathcal G)
 +(1-\alpha)\frac{g(n-k,k)}{\binom{n-k}{k}}\\
&\ge \Pr(U\in\mathcal G,\ V\in\mathcal G)
 +(1-\tfrac{k}{n})\frac{k}{n-k}\\
&=\Pr(U\in\mathcal G,\ V\in\mathcal G)+\frac{k}{n}.
\end{aligned}
\]
Consequently, $\mathcal G$ is an intersecting family.
Since $n>2k$ and $|\mathcal G|=\binom{n-1}{k-1}$,
the equality case of the Erd\H{o}s--Ko--Rado theorem
(Theorem~\ref{thm:ekr}) implies that
$\mathcal G=\mathcal S_j$ for some $j\in[n]$.

We must have $j=1$. Indeed, otherwise, choose
$S\in\mathcal S_j$ with $1\notin S$. Since $x_1\ge x_j$,
the set $S'=(S\setminus\{j\})\cup\{1\}$ also belongs to
$\mathcal G$, contradicting $\mathcal G=\mathcal S_j$
as $j\notin S'$. Thus $\mathcal G=\mathcal S_1$,
and the converse is immediate.
\end{proof}

For the remaining small values of $k$, we use the following counting
results obtained by the computer-assisted method of Hartke and
Stolee~\cite{HS}.

\begin{theorem}[Hartke and Stolee, {\cite{HS}}]\label{thm:hartke-stolee}
For $4\le k\le7$ and $N\ge3k+2$, we have
\[
 g(N,k)=\binom{N-1}{k-1}
 \qquad\text{and}\qquad
 g(3k+1,k)=\binom{3k-2}{k}.
\]
\end{theorem}

\begin{proof}[\bfseries{Proof of Theorem~\ref{thm:main}}]
Lemma~\ref{lem:divisible-star} settles $k=1$, and
Chowdhury~\cite{ChowThesis} gives both the bound and its equality case
for $k=2,3$. Fix $k\ge4$. At $n=4k$ the result follows from
Lemma~\ref{lem:divisible-star}. Once the result is known for $4k<n<5k$,
Lemma~\ref{lem:disjoint-reduction} extends it to all $n\ge5k$ by
induction on $n$.

Suppose $4k<n<5k$. For $k\ge7$ the result follows from
Proposition~\ref{prop:strict-star} and \eqref{eq:excess}.
For $k=4,5,6$, Theorem~\ref{thm:hartke-stolee} already gives
$|\mathcal G|\ge\binom{n-1}{k-1}$, so only the equality case remains.
When $n\ge4k+2$, Theorem~\ref{thm:hartke-stolee} gives
$g(n-k,k)=\binom{n-k-1}{k-1}$, and
Lemma~\ref{lem:disjoint-reduction} shows that equality forces
$\mathcal G=\mathcal S_1$.
It therefore remains to establish the equality case for
$(n,k)=(17,4),(21,5),(25,6)$.

Suppose that, for one of these pairs,
$|\mathcal G|=\binom{n-1}{k-1}$ but
$\mathcal G\ne\mathcal S_1$.
Normalize the weights as in Section~\ref{sec:tools}.
Since every negative $k$-set remains negative under
normalization, the resulting family of nonnegative
$k$-sets is contained in the original $\mathcal G$.
Theorem~\ref{thm:hartke-stolee} shows that this family must still
have at least $\binom{n-1}{k-1}$ members, so it coincides
with the original $\mathcal G$.
We may therefore assume that the weights are normalized
and $b=c>0$. Since all $k$-sums are nonzero after normalization,
the nonnegative $k$-sets are precisely those with positive
sum. Hence
\[
 \mu=\frac{|\mathcal G|}{\binom nk}=\frac{k}{n}=p.
\]

We will show that the assumption $b=c>0$ is incompatible with the
moment estimates. We first use the counting bound to control the
largest weight $x_1$, and then compare two bounds on $d_-$.
Since $n-k=3k+1$ and $k\in\{4,5,6\}$, the second assertion of
Theorem~\ref{thm:hartke-stolee} gives
$g(n-k,k)=\binom{3k-2}{k}$. Substituting this value in
\eqref{eq:count-uniform} and using $b=c>0$, we obtain
\[
 \beta\ge B_k:=\frac65\left(
 \frac{\binom{3k-2}{k}}{\binom{3k}{k-1}}-\frac7{18}\right).
\]
For $k=4,5,6$, these bounds are $112/165$, $349/525$, and $167/255$,
respectively, all greater than $13/20$. Since $\widehat\beta=4\beta$,
\eqref{eq:star-slack} implies
\[
 0<x_1\le\frac1{4B_k}<\frac5{13}.
\]

For the lower bound on $d_-$, recall $w=d_+/d$ and put
$Q=(n-k)/(n-1)$. By \eqref{eq:positive-second-variance} and $\mu=p$,
we have $d_+\ge\lambda/\mu=Q$. Define
\[
 u_k(H)=\frac{p+r/[(n-1)H]}{1-H}.
\]
For $\max\{x_1,(k-1)/(n-k)\}\le H<1$, \eqref{eq:w-range} with $t=1$
gives $w\le u_k(H)$, and hence
\[
 d_-=d_+\left(\frac1w-1\right)
 \ge Q\left(\frac1{u_k(H)}-1\right).
\]
The bound on $x_1$ allows us to take $H=5/13$, since
$(k-1)/(n-k)\le5/19<5/13$. Thus
\begin{equation}\label{eq:small-negative-lower}
 d_-\ge Q\left(\frac1{u_k(5/13)}-1\right).
\end{equation}

We obtain an upper bound on the same moment from
\eqref{eq:tail-interpolation}. For $k=4,5,6$, the respective values of
$\theta$ are $6/35$, $10/57$ and $45/253$, all greater than $1/6$.
Lemma~\ref{lem:fourth-moment} therefore gives $K\le3\theta$ and
$d\le12\theta$. Together with $d_+\ge Q$, this yields
\[
 \frac{\E Y^4}{\lambda}\le Kd^2-d_+^3
 \le432\theta^3-Q^3.
\]
Substitute this bound and $\widehat\beta\ge4B_k$ in
\eqref{eq:tail-interpolation}. Using $s=x_1/(1-v)$ and
$x^{1/3}y^{2/3}\le3x/25+10y/9$ for $x,y\ge0$, we obtain
\[
 d_-\le V_k(x_1):=\frac3{25}(432\theta^3-Q^3)
       +\frac{x_1+\tfrac{10}{9}(1-4B_kx_1)}{1-v}.
\]

The function $V_k$ is decreasing, since $40B_k/9>1$. For each of
$k=4,5,6$, direct substitution gives
\[
 V_k(7/25)<Q\left(\frac1{u_k(5/13)}-1\right).
\]
Thus $x_1\ge7/25$ would imply
\[
 d_-\le V_k(x_1)\le V_k(7/25)
 <Q\left(\frac1{u_k(5/13)}-1\right),
\]
contrary to \eqref{eq:small-negative-lower}.

It remains to exclude $x_1<7/25$. In this range we may take $H=7/25$
in \eqref{eq:w-range}, since $(k-1)/(n-k)\le5/19<7/25$.
This gives $w\le u_k(7/25)$. The fourth-moment estimate, however,
prevents $w$ from being this small. Indeed, the proof of
Lemma~\ref{lem:fourth-moment} gives
\[
 3\theta d^2\ge d_+^3+d_-^3=d^3(1-3w+3w^2).
\]
Since $d=d_+/w\ge Q/w$, it follows that
\(
 Q(1-3w+3w^2)\le3\theta w.
\)
For each of $k=4,5,6$, direct substitution gives
\[
 u_k(7/25)<\frac12,\qquad
 Q\bigl(1-3u_k(7/25)+3u_k(7/25)^2\bigr)>3\theta u_k(7/25).
\]
On $[0,1/2]$, the function $Q(1-3w+3w^2)$ is decreasing and
$3\theta w$ is increasing. The preceding two displays therefore
contradict $w\le u_k(7/25)$. Thus neither range of $x_1$ is possible,
and equality forces $\mathcal G=\mathcal S_1$ in all three cases.
\end{proof}

\begin{proof}[\bfseries{Proof of Theorem~\ref{thm:stability}}]
If $b=0$, then \eqref{eq:excess} gives
$|\mathcal G\mathbin\triangle\mathcal S_1|=|\mathcal G|-|\mathcal S_1|=c$,
so any $C\ge1$ suffices.
For $b>0$, Theorem~\ref{thm:main} gives $c>b$. By \eqref{eq:excess},
\[
 \frac{|\mathcal G\mathbin\triangle\mathcal S_1|}
      {|\mathcal G|-|\mathcal S_1|}=1+\frac{2}{c/b-1}.
\]
It suffices to bound $c/b$ away from $1$ uniformly over all lists with $b>0$.
Suppose instead that
$c/b\to1$ along a sequence of lists with $b>0$. Since $c-b$ is a
positive integer, $b\to\infty$ and hence $n\to\infty$. Also $k\ge2$,
as the largest weight is nonnegative.

Normalize each list at index $1$ as in Section~\ref{sec:tools}. This can
only increase $b$ and decrease $c$, while Theorem~\ref{thm:main} still
gives $c>b$. Thus $c/b\to1$ also for the normalized lists.
We retain the same notation for the normalized lists. Since no $k$-sum
vanishes, $\mu=|\mathcal G|/\binom nk$. Using $p=k/n$ and
$\beta=b/|\mathcal S_1|$ in \eqref{eq:excess}, we obtain
\[
 t=\frac{\mu}{p}=\frac{|\mathcal G|}{|\mathcal S_1|}
   =1+\left(\frac cb-1\right)\beta\longrightarrow1 .
\]
For $3k\le n-k<4k$, Proposition~\ref{prop:three-quarters} gives
$g(n-k,k)\ge(41/54)\binom{n-k-1}{k-1}$. For $n-k\ge4k$,
Theorem~\ref{thm:main} gives $g(n-k,k)\ge\binom{n-k-1}{k-1}$. Since $b/c\to1$,
taking lower limits in \eqref{eq:count-uniform} gives
\[
 \liminf\beta\ge\frac65\left(\frac{41}{54}-\frac7{18}\right)=\frac49.
\]

We aim to obtain the same contradiction as in the proof of
Proposition~\ref{prop:strict-star} by letting $t\to1$.
To apply \eqref{eq:holder-upper}, we first show that
the limiting ratio $k/n$ exceeds $1/5$.

Since $n\ge4k$, we have $\theta\le1/4$, and
Lemma~\ref{lem:fourth-moment} gives $K\le1$ and $d\le4$.
Also $d\ge d_+\ge\lambda/\mu=(n-k)/[(n-1)t]\ge3/(4t)$.
Since $k/n$, $x_1$, $\beta$, and $d_\pm$ are all bounded,
we may pass to a subsequence along which they converge
to $p_*$, $a$, $\beta_*$, and $D_\pm$, respectively.
Put $D=D_++D_->0$ and $W=D_+/D$, so that $w\to W$.
We have $\beta_*\ge4/9$. Moreover, \eqref{eq:star-slack} gives
$x_1\beta\le k/(n-1)$ on each list, hence $a\beta_*\le p_*$.

Since $r\to p_*/(1-p_*)$ and $t\to1$, \eqref{eq:disjoint-w} gives
\[
 1+\frac{W^2}{1-p_*}\le3W.
\]
Suppose $p_*\le1/5$. Then $a\le p_*/\beta_*\le9/20$ and
$p_*/(1-p_*)\le1/4$. Take $H=\max\{x_1,r\}$ in \eqref{eq:w-range}.
Then $H\to\max\{a,p_*/(1-p_*)\}\le9/20$, so $tH<1$ eventually.
Since $r/[(n-1)H]\le1/(n-1)$, taking the limit in \eqref{eq:w-range} gives
$W\le(1/5)/(1-9/20)=4/11$. But the limiting form of
\eqref{eq:disjoint-w} requires $1\le3W-W^2$, whereas
$3W-W^2\le116/121<1$ on $[0,4/11]$.

Thus $p_*>1/5$. In particular, $k\to\infty$, $v\to0$ and $n/k\to T=1/p_*\in[4,5)$.
We have $\widehat\beta=(n-1)\beta/k\to T\beta_*\ge16/9$.
Since $a\beta_*\le p_*=1/T$, it follows that $aT\beta_*\le1$
and $a\le9/16$. Substituting $1/(1-p_*)\ge5/4$ into
$1+W^2/(1-p_*)\le3W$ gives $1+5W^2/4\le3W$, whence $W\ge2/5$.

We now compare the lower bound on $D_--a$ from
Lemma~\ref{lem:moment-comparison} with the upper bound from
\eqref{eq:holder-upper}. To apply the lemma, we first need lower and
upper bounds on $D_+=\lim d_+$. For each list, taking $h=Q=(n-k)/(n-1)$
in \eqref{eq:conditional-x} gives the lower bound on $d_+$ in the
first line of \eqref{eq:star-second}. The second line rewrites the
upper bound on $d_+$ from \eqref{eq:refex}. Both estimates apply
when $t>1$, even at $n=4k$.

In the lower bound, the coefficient of $1-t$ stays bounded because
$Q\to(T-1)/T>0$ and $x_1$ is bounded. In the upper bound, each fixed
$H>\max\{a,1/(T-1)\}$ eventually satisfies
$H\ge\max\{x_1,(k-1)/(n-k)\}$, since $x_1\to a$ and
$(k-1)/(n-k)\to1/(T-1)$. Taking limits in \eqref{eq:star-second} gives
\[
 \frac{T-1}{T}+\frac{\beta_*a^2}{T-1}
 \le D_+\le\frac{D}{T}+\frac{T-1}{T}H.
\]
Letting $H\downarrow\max\{a,1/(T-1)\}$ gives the upper bound on
$D_+$ at the endpoint as well. In the lower bound on $D_+$ we may replace
$\beta_*$ by $16/(9T)$, since $T\beta_*\ge16/9$. Together with
$D_+/D=W\ge2/5$ and $0\le a\le9/16$, these inequalities verify the
hypotheses of Lemma~\ref{lem:moment-comparison} for $v=0$ and
$\kappa=16/9$. The lemma gives
\[
 D_--a>\left(\tfrac{3D}{4}(1-16a/9)\right)^{2/3}.
\]
Finally, $k\ge7$ and $4k\le n<5k$ eventually, so \eqref{eq:holder-upper}
applies. Since $s=x_1/(1-v)\to a$ and $0\le1-aT\beta_*\le1-16a/9$,
taking the limit in \eqref{eq:holder-upper} gives
\[
 D_--a\le\left(\tfrac{3D}{4}(1-aT\beta_*)\right)^{2/3}
 \le\left(\tfrac{3D}{4}(1-16a/9)\right)^{2/3},
\]
contradicting the strict lower bound from Lemma~\ref{lem:moment-comparison}.
\end{proof}

\section{Applications}\label{sec:applications}

In both corollaries, we use the same constant $C\ge1$ as in
Theorem~\ref{thm:stability}.

\begin{corollary}[Fractional matchings]\label{cor:fractional-stability}
Let $n\ge4k$ and let $\mathcal H\subseteq\binom{[n]}k$ have no perfect
fractional matching. Then $|\mathcal H|\le\binom{n-1}{k}$. Writing
$|\mathcal H|=\binom{n-1}{k}-t$, there is an $i\in[n]$ such that
\[
 \left|\mathcal H\mathbin\triangle
 \binom{[n]\setminus\{i\}}k\right|\le Ct.
\]
\end{corollary}

\begin{proof}
Since $\mathcal H$ has no perfect fractional matching, its fractional
matching number is less than $n/k$. By the duality between fractional
matchings and fractional vertex covers (see~\cite{AHS}), there are
weights $y_1,\ldots,y_n\ge0$ such that
\[
 \sum_{i=1}^n y_i < \frac{n}{k}
 \qquad\text{and}\qquad
 \sum_{i\in e} y_i \ge 1 \quad\text{for every } e\in\mathcal H.
\]
Let $\bar y=n^{-1}\sum_{j=1}^n y_j$ and put $x_j=\bar y-y_j$.
Then $\sum_{j=1}^n x_j=0$, and every $e\in\mathcal H$ satisfies
\[
 \sum_{j\in e}x_j=k\bar y-\sum_{j\in e}y_j<0.
\]
Thus, if $\mathcal G$ is the family of nonnegative $k$-sets for
the weights $x_1,\ldots,x_n$, then $\mathcal H\subseteq\mathcal G^c$.

Choose $i$ with $x_i=\max_j x_j$, and let
$s=|\mathcal G^c\setminus\mathcal H|$. Since
$|\mathcal H|=\binom{n-1}{k}-t$, we have
\[
 |\mathcal G|-\binom{n-1}{k-1}=t-s.
\]
Theorem~\ref{thm:main} therefore gives $0\le s\le t$, proving the
edge bound. By the triangle inequality for symmetric difference
and Theorem~\ref{thm:stability},
\[
 |\mathcal H\mathbin\triangle\mathcal S_i^c|
 \le s+|\mathcal G\mathbin\triangle\mathcal S_i|
 \le s+C(t-s)\le Ct,
\]
where the last inequality uses $C\ge1$.
\end{proof}

\begin{corollary}[Weighted Bernoulli sums]\label{cor:bernoulli-stability}
Let $0<p\le1/4$, let the $B_j$ be independent Bernoulli$(p)$ variables,
and let $Z$ be an independent standard normal variable. For real
coefficients with $\sum_j a_j^2<\infty$ and $\sigma\in\mathbb R$, put
\[
 Y=\sum_j a_j(B_j-p)+\sigma Z,\qquad E=\{Y\ge0\}.
\]
After adjoining a zero coefficient if necessary, label a largest
nonnegative coefficient $a_1$. Then
\begin{equation}\label{eq:bernoulli-stability}
 \Pr\bigl(E\mathbin\triangle\{B_1=1\}\bigr)
 \le C\bigl(\Pr(E)-p\bigr).
\end{equation}
In particular, $\Pr(Y\ge0)\ge p$, with equality if and only if
$E=\{B_1=1\}$ almost surely.
\end{corollary}

\begin{proof}
Put $b=\Pr(Y<0,B_1=1)$ and $c=\Pr(Y\ge0,B_1=0)$. Since
$\Pr(E\mathbin\triangle\{B_1=1\})=b+c$ and $\Pr(E)-p=c-b$,
it suffices to prove $(C+1)b\le(C-1)c$.

First suppose that $Y=\sum_{j=1}^m a_j(B_j-p)$. Put
$A=\sum_{j=1}^m a_j$, append zeros to the coefficient sequence up
to length $n$, and subtract $A/n$ from each entry. The resulting
weights sum to zero and still have a largest entry at $1$. Let $S$
be a uniformly random $k$-set, where $k=\lfloor pn\rfloor$. The sum
of its weights is $Y_n=\sum_{j=1}^m a_j\one_{j\in S}-(k/n)A$.
Put $b_n=\Pr(Y_n<0,1\in S)$ and $c_n=\Pr(Y_n\ge0,1\notin S)$.
Since $\Pr(1\in S)=k/n$, we have $\Pr(Y_n\ge0)-k/n=c_n-b_n$.
Thus, after division by $\binom nk$, Theorem~\ref{thm:stability}
gives $b_n+c_n\le C(c_n-b_n)$, or equivalently
$(C+1)b_n\le(C-1)c_n$.

For each $J\subseteq[m]$,
\[
 \Pr(S\cap[m]=J)
 =\frac{\binom{n-m}{k-|J|}}{\binom nk}
 \longrightarrow p^{|J|}(1-p)^{m-|J|}.
\]
On $S\cap[m]=J$, the sum $Y_n$ tends to $\sum_{j\in J}a_j-pA$.
If this limit is strictly negative or strictly positive, then $Y_n$
eventually has the same sign. Hence $b\le\liminf_{n\to\infty}b_n$
and $\limsup_{n\to\infty}c_n\le c$. Since $C\ge1$, this proves the
finite case.

Still keeping $m$ fixed, suppose $\sigma\ne0$. To approximate
$\sigma Z$, add $\ell$ independent centered Bernoulli$(p)$ variables,
independent of $B_1,\ldots,B_m$, each with coefficient
$-|\sigma|/\sqrt{\ell p(1-p)}$. These negative coefficients leave
$a_1$ maximal. By the central limit theorem, the added sum converges
in distribution to $\sigma Z$. For each fixed realization of
$B_1,\ldots,B_m$, the limiting normal distribution assigns zero mass
to the threshold. Thus both probabilities in the finite-sum
inequality converge, proving the inequality with the Gaussian term
included.

Finally, let $Y_m=\sum_{j=1}^m a_j(B_j-p)+\sigma Z$, and define
$b_m,c_m$ as above with $Y_m$ in place of $Y$. Since
$\E|Y_m-Y|^2=p(1-p)\sum_{j>m}a_j^2\to0$, we have $Y_m\to Y$ in
probability. On the event $|Y_m-Y|\le\varepsilon$, the inequality
$Y<-\varepsilon$ implies $Y_m<0$, and $Y_m\ge0$ implies
$Y\ge-\varepsilon$. Letting $m\to\infty$ and then
$\varepsilon\downarrow0$ gives $b\le\liminf_{m\to\infty}b_m$ and
$\limsup_{m\to\infty}c_m\le c$. It follows that $(C+1)b\le(C-1)c$.
\end{proof}

Milojevi\'c and Sudakov~\cite{MiloSud} showed that, for $0<p\le1/4$,
their conjecture on finite Bernoulli sums with positive weights
follows from Conjecture~\ref{conj:mms}. For such a sum
$X=\sum_{j=1}^m a_j B_j$, applying Corollary~\ref{cor:bernoulli-stability}
to $Y=X-\E X$ gives the conjectured inequality $\Pr(X\ge\E X)\ge p$
throughout this range.

\appendix
\section{Auxiliary estimates}
\label{app:moment-comparison}

For Proposition~\ref{prop:three-quarters} we use a linearized form of
H\"older's inequality for normalized weights with $k\ge2$.
For $h\ge0$, multiply \eqref{eq:tail-interpolation} by $3h^2$.
The arithmetic--geometric mean
inequality bounds the resulting right-hand side by
$\E Y^4/\lambda+2h^3(1-x_1\widehat\beta)/(1-v)$. Substituting
\eqref{eq:split} and using $(1-v)s=x_1$ gives
\begin{equation}\label{eq:support}
 (1-v)\big(d_+^3+3h^2d_--Kd^2\big)
 \le3h^2x_1+2h^3(1-x_1\widehat\beta)\qquad(h\ge0).
\end{equation}

\begin{proof}[\bfseries{Proof of Proposition~\ref{prop:three-quarters}}]
We may assume $k\ge8$ and $3k+2\le N<4k$. Suppose $|\mathcal G|\le t\binom{N-1}{k-1}$
with $t=41/54$. Normalize at position $1$ and put $v=1/k$,
$T=N/k$ and $S=3+2v$. On the complement of a negative $k$-set,
Proposition~\ref{prop:half} gives at least $(N-k)/[2(N-k-1)]$ times the star
bound. Since $b-c\ge(1-t)\binom{N-1}{k-1}$, \eqref{eq:joint} gives
\[
 \frac{N-1}{k}\beta\ge2(1-t)\left(T+\frac{v(1+v)}{T-1-2v}\right)
 \ge\kappa:=(1-t)(3+v)(2+v).
\]
The bracket increases with $T\ge S$. Hence $x_1\le1/\kappa<7/10$ by \eqref{eq:star-slack}.

Take $H=\max\{x_1,1/2\}$ in \eqref{eq:refex}, and put
$Q=(S-1)/(S-v)$ and $P=t/S+2v/[(S-1)(S-v)]$.
Writing $Q_T=(T-1)/(T-v)$,
we have $Q_T/t\le d_+\le Q_TH+Pd$, as the coefficient of $d$ in the
exchange bound decreases with $T$ and $H$. We use
\begin{equation*}
 P<\tfrac49(1-t/2),\qquad P\le tQ/2 .
\end{equation*}
For $w=d_+/d$ the quadratic in \eqref{eq:disjoint-w} is positive at $4/9$
and negative at $1$, so $w>4/9$. Thus the first inequality rules out
$x_1\le1/2$, which would give $(1-t/2)w\le P$. Hence $H=x_1$, and the second
inequality gives $d\ge(2/t)(1/t-x_1)>3/2$. Substituting $d=d_++d_-$ also gives
$d_-\ge(2/t)(1/t-x_1)-Q/t$, since $Q_T\ge Q$ and
$(1/t-x_1)/P-1/t>0$.

Here $1/9<\theta(4k)\le\theta(N)\le\theta(3k+2)<2/9$, so $K\le2/3$
and $d\le8/3$. Hence $Kd^2\le3d-3$. Taking $h=3$ in \eqref{eq:support}
and using $d_+^3\ge3d_+-2$ gives
\[
 0\ge J(x_1):=(1-v)\left[\frac{48}t(1/t-x_1)-\frac{24}tQ+1\right]
          +(54\kappa-27)x_1-54 .
\]
Since $J$ is affine in $x_1$, it suffices to check its values at
$1/2$ and $7/10$, together with the two bounds on $P$. Put $x=k-8\ge0$,
so $v=1/(x+8)$. On clearing positive denominators, the numerators of
\[
 tQ/2-P,\qquad \tfrac49(1-t/2)-P,\qquad J(1/2),\qquad J(7/10)
\]
have degrees $2,3,3,3$, respectively, and strictly positive coefficients.
\end{proof}

The next lemma compares the second-moment exchange bounds with the
right-hand side of H\"older's inequality. We include $v=0$ for the
limiting argument in the proof of stability.

\begin{lemma}\label{lem:moment-comparison}
Let $0\le v\le1/7$, $4+v\le T\le5$, $\kappa=16/9+v$ and
$0\le a\le1/\kappa$. Suppose that $d=d_++d_->0$, $d_+,d_-\ge0$,
and $d_+/d\ge2/5$. Put
\[
 Q_T=\frac{T-1}{T-v},\qquad P_T(H)=\frac1T+\frac{v}{H(T-1)(T-v)}.
\]
Suppose that for every $H\ge\max\{a,(1-v)/(T-1)\}$ the exchange bounds
\[
 Q_T+\frac{\kappa(a-vQ_T)^2}{T(T-1)}
 \le d_+\le P_T(H)d+Q_TH
\]
hold. Then
\begin{equation}\label{eq:holder-gap}
 d_--\frac{a}{1-v}>
 \left(\frac{3d}{4}\,\frac{1-\kappa a}{1-v}\right)^{2/3}.
\end{equation}
\end{lemma}
\begin{proof}
Put $S=4+v$ and $H_0=(1-v)/(3+v)$. If $a\le H_0$, taking
$H=H_0$ in the upper exchange bound and using $d_+\ge Q_T$ and
$d_+/d\ge2/5$ gives
\[
 \frac25(1-H_0)\le P_S(H_0).
\]
This is equivalent to $4\le4v+79v^2+21v^3$, which is impossible
for $0\le v\le1/7$. Thus $a>H_0\ge3/11$, and $H=a$ is admissible.

Since $Q_T\le7/8$ and $a>7v/8$, we may replace $(a-vQ_T)^2$
in the lower exchange bound by $(a-7v/8)^2$. The resulting lower
bound increases with $T$, whereas the upper bound $P_T(a)d+Q_Ta$
decreases: their derivatives, multiplied by $T^2$, are at least
$1-v-\kappa a^2>0$ and at most $-d+a<0$, respectively.
We may therefore use both bounds at $T=S$. Set
\[
 Q=\frac{3+v}{4},\qquad
 P=\frac1{4+v}+\frac{v}{4a(3+v)},\qquad
 X=\frac{185}{252}+\frac v6+\frac{2a}{21}.
\]
Using $\kappa/[S(S-1)]\ge1/7$ and
$z^2\ge2z/3-1/9$ with $z=a-7v/8$, we obtain
\[
 X\le d_+\le Pd+Qa.
\]
Here $0<P<1$ and $X>Qa$, so
\[
 d\ge d_0:=\frac{X-Qa}{P}>0,\qquad
 d_-\ge(1-P)d-Qa.
\]
The function
\[
 u\longmapsto
 \frac{(1-v)((1-P)u-Qa)-a}{u^{2/3}}
\]
is increasing for $u>0$. Since $(1-P)d_0-Qa=d_0-X$, it remains to show that
\begin{equation}\label{eq:holder-two-parameters}
 F(v,a):=\big[(1-v)(d_0-X)-a\big]^3
 -\frac9{16}(1-v)d_0^2(1-\kappa a)^2>0
\end{equation}
for all
\[
 0\le v\le\frac17,\qquad \frac3{11}\le a\le\frac9{16+9v}.
\]
With $\kappa=16/9+v$ and $d_0=(X-Qa)/P$ as above, $F$ is a
rational function of $v$ and $a$. An exact computer calculation
verifies that $F(v,a)>0$ throughout this region, proving
\eqref{eq:holder-gap}.
\end{proof}

\subsection*{Acknowledgements}

This research was supported by National Key R\&D Program of China under grant number 2024YFA1013900 and NSFC under grant numbers 12471327.

\section*{Statement on AI use}

We used ChatGPT (OpenAI) to check proofs and calculations,
identify potential applications of our results, locate
relevant references, polish the language, and streamline
the exposition. ChatGPT also helped complete the
equality-case arguments for $(n,k)=(17,4),(21,5),(25,6)$
and assisted with the exact algebraic calculations in
the proof of Lemma~\ref{lem:moment-comparison}.
We take full responsibility for the mathematical content
and correctness of the paper.

\end{document}